\documentclass[10pt,twocolumn,twoside]{IEEEtran}

\usepackage{amssymb,amsmath,amsthm}
\usepackage{bm,mathrsfs}
\usepackage{enumerate}
\usepackage{hyperref}
\usepackage{cleveref}
\usepackage{graphicx,tikz}
\usepackage{algorithm,algorithmic}
\usepackage[tight,footnotesize]{subfigure}
\usepackage{multirow}
\usepackage{rotating}
\usepackage{dsfont}
\usepackage{booktabs}
\usepackage{color}

\DeclareMathOperator*{\argmin}{arg\,min}

\providecommand{\norm}[1]{\left\rVert#1\right\rVert}
\providecommand{\set}[1]{ \left\{ #1  \right\}  }

\providecommand{\parenth}[1]{\left( #1 \right) }

\begin{document}
%--------------------------------------------------------------

\title{Discretization and Optimization using Graphs: \\
 One-Dimensional Algorithm }
\author{
John Paul Ward

\thanks{The author is a member of the  Department of Mathematics and Statistics at North Carolina A\&T State University.}
\thanks{Python code is provided on arXiv.}

}

\maketitle

\ifCLASSOPTIONpeerreview
\begin{center}
 \vspace{-0.6cm}
 Authors' contact information:\\
 \vspace{0.2cm}
 North Carolina A\& T State University\\
 Mathematics Department\\
 Greensboro, North Carolina.
 \vspace{0.2cm}
\end{center}
\fi

\begin{abstract}
\boldmath
We consider the problem of discretizing one-dimensional, real-valued functions as graphs.
The goal is to find a small set of points, from which we can approximate the remaining function values. 
The method for approximating the unknown values is interpolation on a discrete graph structure.
From the discrete graph structure, we build a refined approximation to the function over its domain.
This fine approximation can then be used for problems such as optimization, which we illustrate by identifying local minima.

\end{abstract}

\ifCLASSOPTIONjournal
\fi

\ifCLASSOPTIONpeerreview
 \begin{center} \bfseries EDICS Category: SAS-MALN,  MLAS-GRAPH \end{center}
\fi

\IEEEpeerreviewmaketitle

%----------------------------------------------------------------------------------------
%----------------------------------------------------------------------------------------

\section{Introduction}
\label{sec:introduction}

Graphs provide a versatile discrete structure for approximation on a variety of continuous domains. They are also useful for creating a structured representation of data that may not have a natural, underlying continuous domain. Here, we focus on the former and investigate how well graphs can approximate functions on bounded, one-dimensional domains.

Our primary goal is to create discrete representations of functions. We are motivated by the problem of optimizing functions that are difficult to compute (or approximate) even pointwise \cite{ampgo}.  By creating a discrete representation, we can then use spline interpolation to fill out approximate function values at any point of the domain.  

To illustrate our approach, we consider sample functions from \cite{ampgo} that have a wide range in smoothness.
For the experiments, we run our algorithm with a generic set of parameters.  The algorithm adapts to each function and provides a discrete representation. From this, we build a fine-grid approximation over the domain and compare with the ground truth. We measure the average $L_2$ error and compare the global minimum of each function with the list of local minima in the approximation.

Interpolation and approximation utilizing graph domains has been studied by several authors \cite{cavoretto2021partition,emirov2022polynomial,erb2022graph,jiang2019nonsubsampled,pesenson00,pesenson08,pesenson09,pesenson10s,pesenson10r,pesenson2021graph,smola03,ward2020interpolating}.
These works cover both the general theory of approximation as well as specific applications utilizing graph structures. 
Additionally, for a general overview of signal processing on graphs, including some specific applications, see 
\cite{ortega2018graph,ortega2022introduction,ricaud2019fourier,shuman13}.

While this work focuses on approximation of functions on one-dimensional continuous domains, we point out that the generality of the graph structure allows for extensions in several directions. Our main goal is to establish a starting point for further investigation and begin the development of an extensible algorithm that is implemented in Python.  
We lay out the basic graph structure under consideration in \Cref{sec:background}. 
The discretization algorithm and experiments are discussed in \Cref{sec:discretization}. 
Optimization is discussed in \Cref{sec:optimization}.
Finally, we conclude in \Cref{sec:summary}.

\subsection{Background and preliminaries}
\label{sec:background}
We consider the problem of discretizing bounded functions on bounded domains, i.e. $f:[a,b]\to \mathbb{R}$. The discrete form of the function is a graph
\begin{equation}
\mathcal{G} = \set{\mathcal{V},E\subset \mathcal{V}\times \mathcal{V},w}
\end{equation}
The vertex set $\mathcal{V}$ contains at least two nodes. These are points in our domain $[a,b]$. The edge set is $E$, which connect sequential points from the vertex set. We assume points are not repeated in the vertex set. The weight function $w:E\rightarrow \mathbb{R}_{>0}$ specifies the closeness of two vertices. For adjacent vertices (adjacent points) $a\leq v < w \leq b$ the weight is 
\begin{equation}
   \frac{1}{\text{dist}(v,w)} =  \frac{1}{w-v}
\end{equation}

The adjacency matrix $A$ of the graph is an $N\times N$ matrix, where $N$ is the number of vertices in the graph. The columns and rows are indexed by the vertices. The $i,j$ entry is zero for pairs of vertices that are not connected by an edge. If vertex $i$ and vertex $j$ are connected by an edge, the $i,j$ entry is the weight of the edge connecting them.

The degree matrix $D$ is a diagonal matrix whose entries are the row sums of the adjacency matrix. The Laplacian on the graph is the matrix $L=D-A$.

\subsection{Interpolation on graphs}

Interpolation on the graph is defined as follows. 
We assume there are known and unknown function values.
We sort the vertices so that the locations on the unknown vertices appear first.
The Laplacian is partitioned as $L = (L_u, L_k)$, where $L_u$ is the submatrix of the Laplacian corresponding to the unknown vertices, 
and $L_k$ is the submatrix of the Laplacian corresponding to the known vertices.  
The function values are represented in a vector 
\begin{equation}
F 
= 
\parenth{
\begin{matrix}
F_u \\
F_k
\end{matrix}
}
\end{equation}
where $F_k$ are the known function values, and $F_u$ is the vector of function values to be computed/approximated.

The spline approximation is found by solving the optimization problem
\begin{align}
\argmin_{F_u}
\norm{(L_u,L_k) 
\parenth{
\begin{matrix}
F_u \\
F_k
\end{matrix}
}
}_2
&=
\argmin_{F_u}
\norm{L_uF_u +L_kF_k}_2
\end{align}
This is done using a standard least-squares solver. For the experiments presented here, we use the NumPy solver \cite{harris2020}.

%----------------------------------------------------------------------------------------
%----------------------------------------------------------------------------------------

\section{Discretization}
\label{sec:discretization}

\subsection{Algorithm}
\label{sec:algorithm}

There are essentially two parts to the discretization algorithm. First, we build a set of sample points from our function. Second, we reconstruct the function on a fine grid on the domain, using only the points evaluated in the first part. The second part is a straightforward interpolation problem from known values. The challenging part is identifying the \emph{best} points in part one so that we end up with a good approximation to the function on its domain.

We assume that the function is computationally expensive to evaluate, so we aim to perform a minimal number of function evaluations in part one, while still having good accuracy in part two. We also assume that we have no knowledge of the function (smoothness etc.) except its domain and the points that we evaluate. 

The algorithm for part one can be summarized as follows:
\begin{enumerate}[a)]
\item
Evaluate the function on a sparse grid spanning the domain
\item
Define a refined grid, and compare the true function values to the approximation resulting from interpolation
\item
Note the locations where the interpolated value is outside a defined tolerance of the true function value
\item
Repeat while there are points where interpolation fails:
\begin{itemize}
    \item 
    Refine the grid in a neighborhood where interpolation fails
    \item
    Compare the true function values to the approximation resulting from interpolation
    \item
    Note the locations where the interpolated value is outside a defined tolerance of the true function value
\end{itemize}

\end{enumerate}

This process returns an array of true function values on a non-uniform grid.  As with many interpolation algorithms, iterpolation is more likely to fail either near the boundary of the domain or at locations where the function is less smooth. Our algorithm uses this to decide which points need to be evaluated.

%----------------------------------------------------------------------------------------
%----------------------------------------------------------------------------------------

\subsection{Experiments -- Discretization}
\label{sec:experiment}

We illustrate our algorithm with examples of both smooth and nonsmooth functions. These examples com from \cite{ampgo}. The functions and domains are shown in \Cref{tab:functions}.

Our implementation is in Python. It takes the following inputs:
\begin{itemize}
\item 
function to be discretized
\item
function domain (bounded)
\item
\emph{data} -- initial grid for evaluation
\item
\emph{ref\_max} -- maximum number of refinements
\item
\emph{tol} -- tolerance for successful pointwise interpolation
\end{itemize}
The following are returned in a dictionary:
 \begin{itemize}
\item 
\emph{nfev} total number of computed function values
\item
discrete approximation -- array of computed values with domain locations 
\end{itemize}

\renewcommand{\arraystretch}{1.2}
\begin{table}[ht!]
\centering
\begin{tabular}{c||c|c}
index	  & domain & function   \\ \hline \hline
2 
& $[2.7,7.5]$ 
& ${\displaystyle  \sin(x)+\sin\parenth{\frac{10x}{3}}   }$ \\ \hline 
3
& $[-10,10]$ 
& ${\displaystyle  -\sum_{k=1}^5 k\sin\parenth{ (k+1)x+k }   }$ \\ \hline 
4
& $[1.9,3.9]$ 
& ${\displaystyle  -(16x^2-24x+5)e^{-x}   }$ \\ 
\hline 
5 
& $[0,1.2]$ 
& ${\displaystyle  -(1.4-3x)\sin(18x)  }$ \\ \hline 
6
& $[-10,10]$ 
& ${\displaystyle -\parenth{x+\sin(x)}e^{-x^2}   }$ \\ \hline 
7
& $[2.7,7.5]$ 
& ${\displaystyle  \sin(x)+\sin\parenth{\frac{10x}{3}}+\log(x)-0.84x+3  }$ \\ 
\hline 
8
& $[-10,10]$ 
& ${\displaystyle  -\sum_{k=1}^5 k\cos\parenth{ (k+1)x+k }   }$ \\ \hline 
9
& $[3.1,20.4]$ 
& ${\displaystyle  \sin(x)+\sin\parenth{\frac{2x}{3}}   }$ \\ \hline 
10
& $[0,10]$ 
& ${\displaystyle  -x\sin(x)   }$ \\ 
\hline 
11
& $\left[-\frac{\pi}{2},2\pi \right]$ 
& ${\displaystyle  2\cos(x)+\cos(2x)   }$ \\ 
\hline 
12
& $\left[0,2\pi \right]$ 
& ${\displaystyle  \sin^3(x)+\cos^3(x)   }$ \\ 
\hline 
13
& $\left[0.001,0.99 \right]$ 
& ${\displaystyle  -x^{2/3}-(1-x^2)^{1/3}  }$ \\ 
\hline 
14
& $\left[0,4 \right]$ 
& ${\displaystyle  -e^{-x}\sin(2\pi x)  }$ \\ 
\hline 
15
& $\left[-5,5 \right]$ 
& ${\displaystyle  \frac{x^2-5x+6}{x^2+1} }$ \\ 
\hline
18
& $\left[0,6 \right]$ 
& ${\displaystyle  
\begin{cases}
(x-2)^2, & x\leq 3 \\
2\log(x-2)+1, & \text{otherwise}
\end{cases}
}$ \\ 
\hline
20
& $\left[-10,10 \right]$ 
& ${\displaystyle -(x-\sin(x))e^{-x^2} }$ \\ 
\hline
21
& $\left[0,10 \right]$ 
& ${\displaystyle x\sin(x)+x\cos(2x) }$ \\ 
\hline
22
& $\left[0,20 \right]$ 
& ${\displaystyle e^{-3x}-\sin^3(x) }$ \\ 
\hline
\end{tabular}
\caption{Functions used in the experiments \cite{ampgo}.}
\label{tab:functions}
\end{table}

\begin{table}[ht!]
\centering
\begin{tabular}{c||c|c|}
index	  & nfev & Average $L_2$ error    \\ \hline \hline
2 & 35 & 1.400e-03  \\ \hline 
3 & 65 & 2.434e-02  \\ \hline 
4 & 23 & 6.290e-04  \\ \hline 
5 & 39 & 6.660e-03  \\ \hline 
6 & 29 & 9.393e-05  \\ \hline 
7 & 35 & 1.395e-03  \\ \hline 
8 & 65 & 2.394e-02  \\ \hline 
9 & 33 & 4.471e-04  \\ \hline 
10 & 25 & 2.668e-03  \\ \hline 
11 & 23 & 7.849e-04  \\ \hline 
12 & 33 & 2.968e-04  \\ \hline 
13 & 31 & 1.424e-03  \\ \hline 
14 & 41 & 8.279e-04  \\ \hline 
15 & 31 & 4.590e-04  \\ \hline 
18 & 27 & 7.649e-04  \\ \hline 
20 & 37 & 6.548e-06  \\ \hline 
21 & 35 & 4.211e-03  \\ \hline 
22 & 63 & 5.757e-04  \\ \hline 
\end{tabular}
\caption{Reported number of function evaluations(nfev) and average $L_2$ error for interpolated approximation for the stated list of functions. 
}
\label{tab:disc_results_generic}
\end{table}

\begin{figure}[!htb]
\centering
     \includegraphics[scale=0.4]{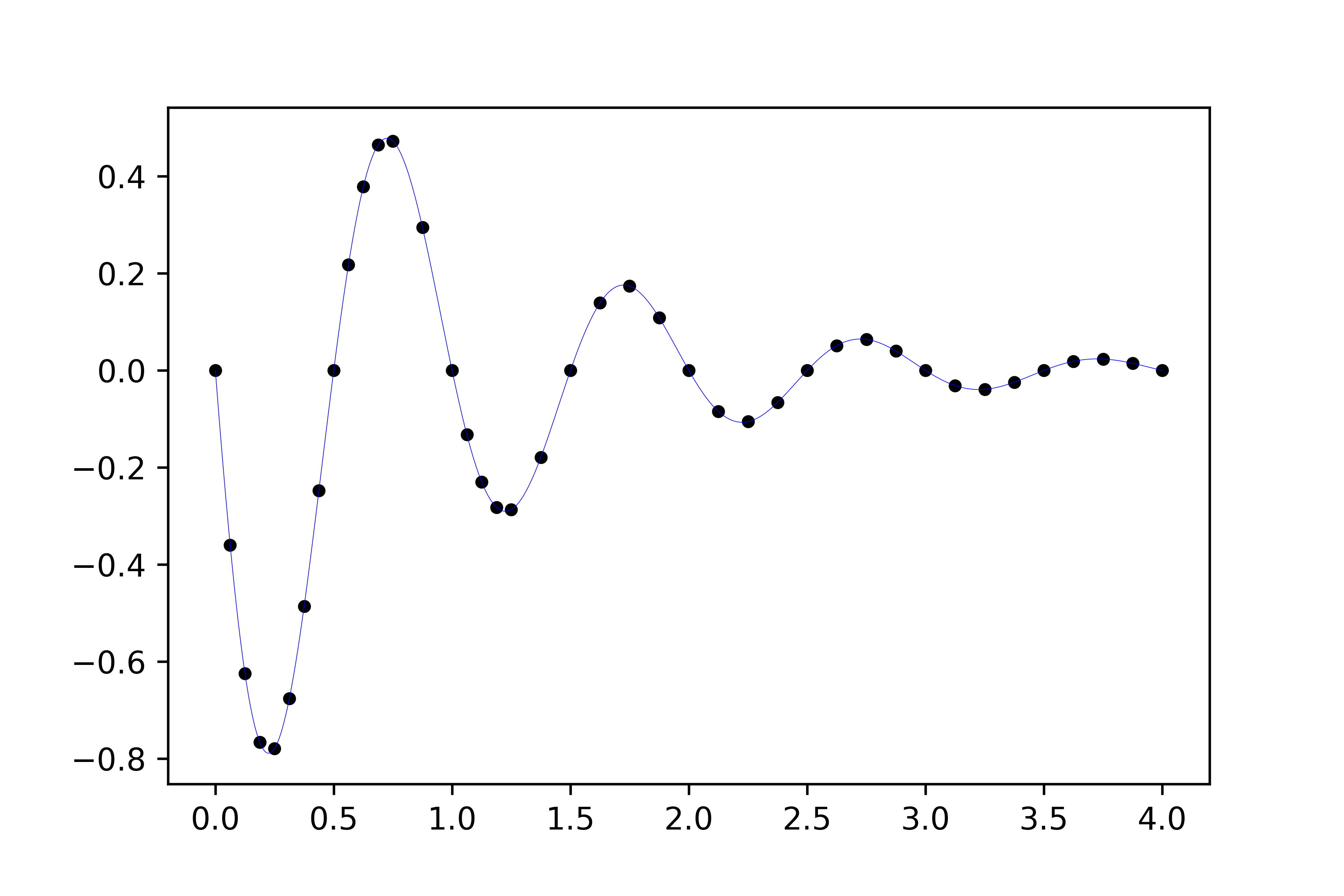}
     
     \includegraphics[scale=0.4]{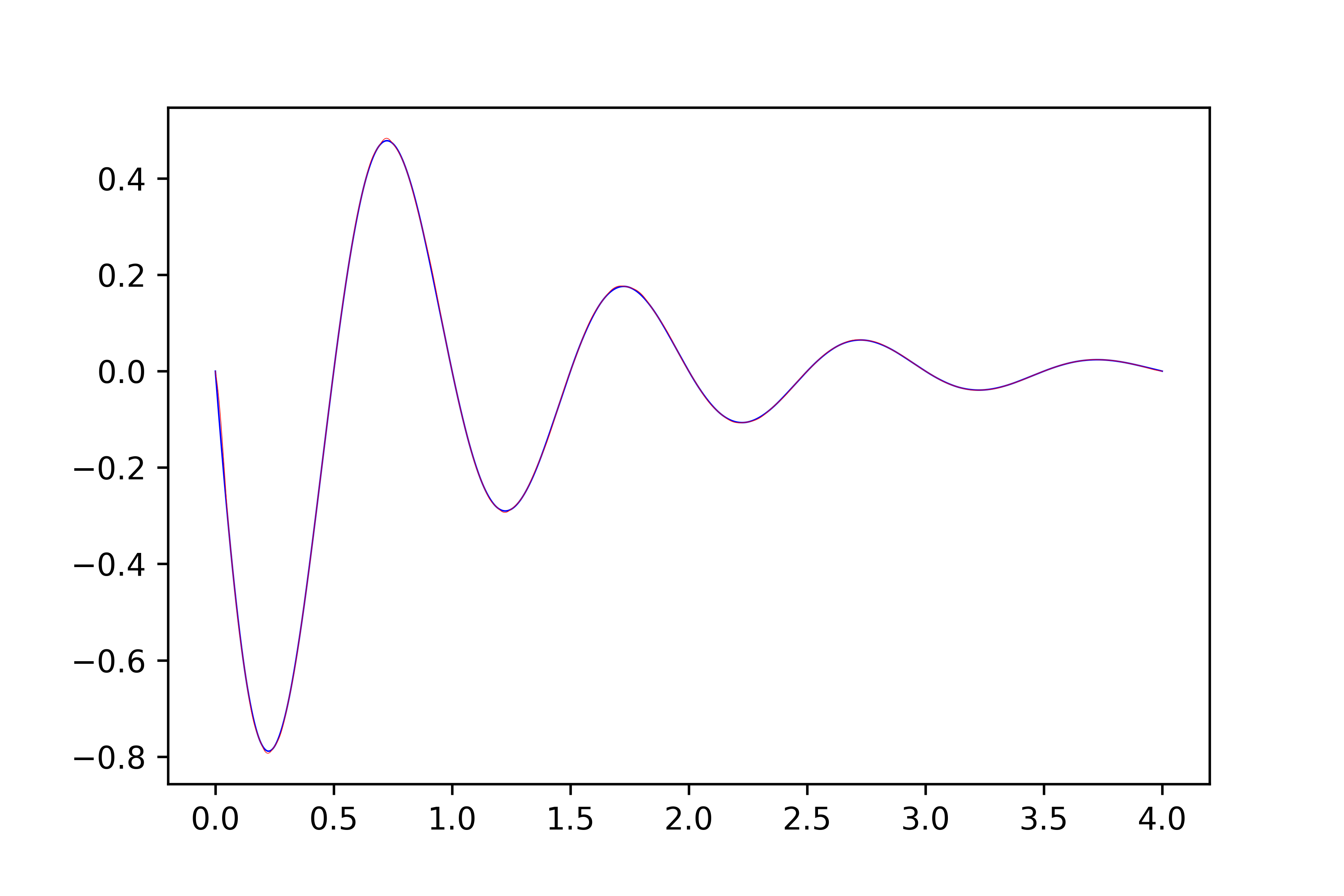}

     \includegraphics[scale=0.4]{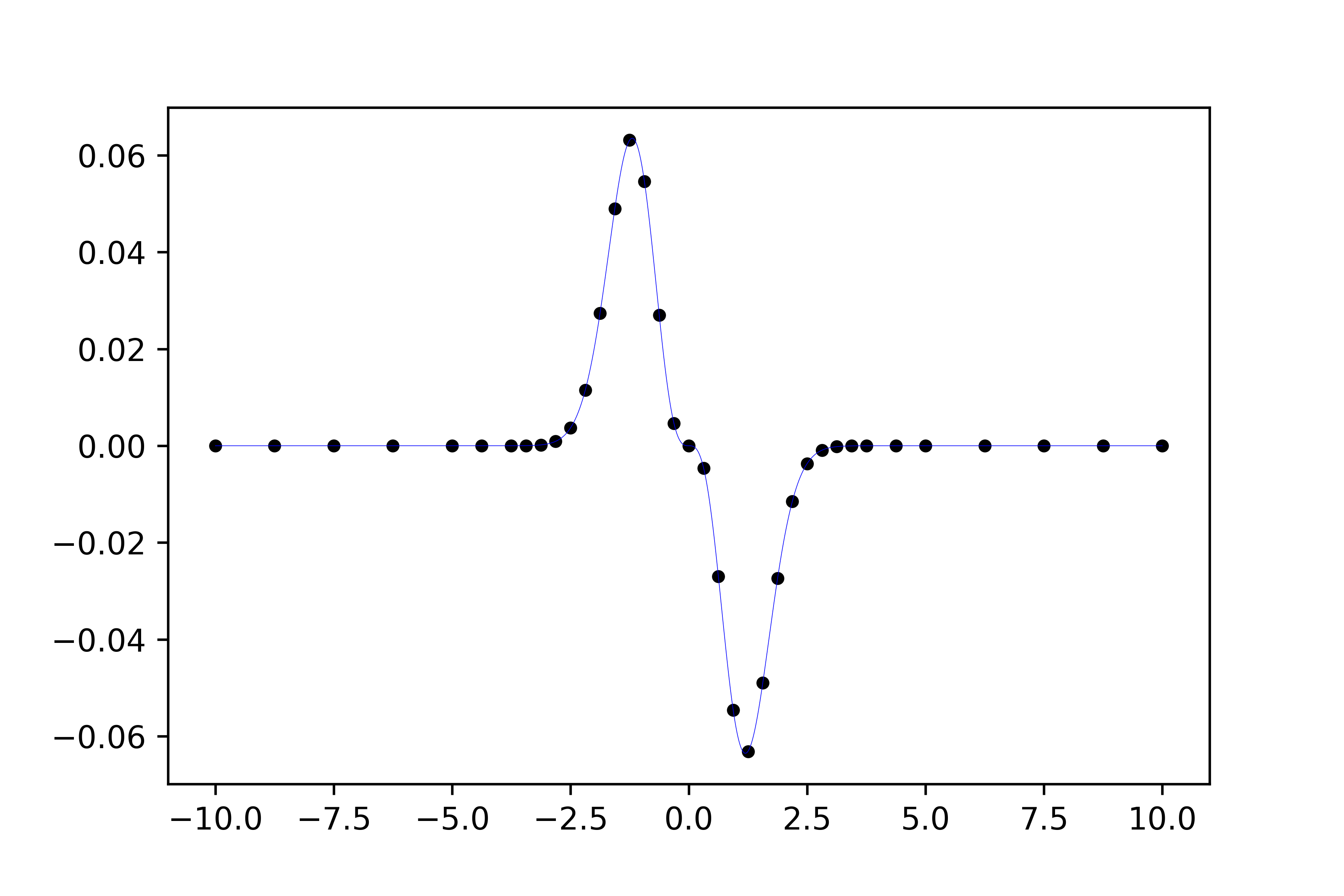}

     \includegraphics[scale=0.4]{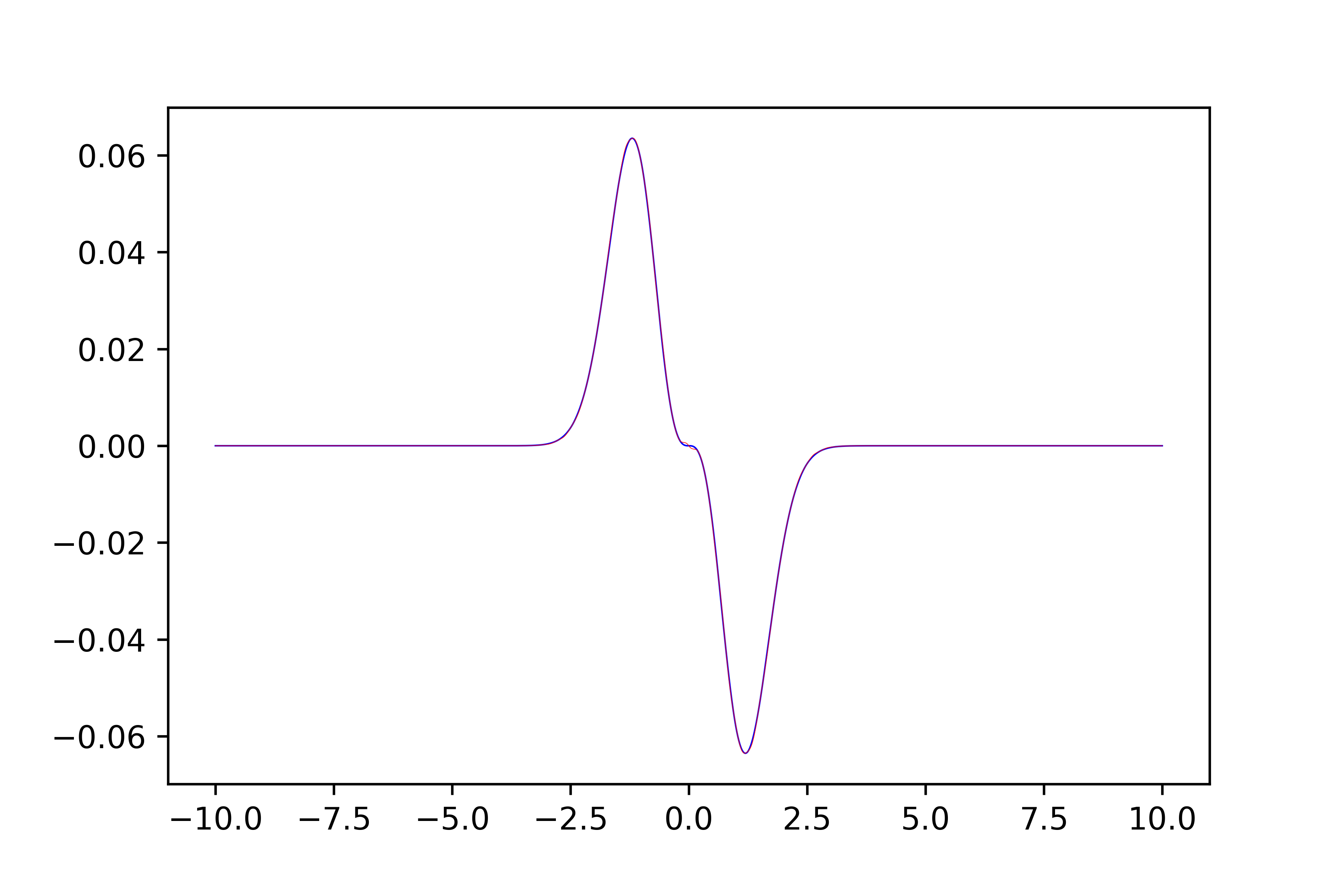}
\caption{
Function evaluation locations and fully interpolated function values for two specific functions, those with indices 14 and 20.}
\label{fig:plots}
\end{figure}

In the experiments, we use an initial default grid of  $9$ uniformly spaced points, $ref\_max=3$, and the point-wise tolerance is $10^{-2}$.
The results in \Cref{tab:disc_results_generic} show that we achieve a good approximation with a relatively small number of initial points and only a few refinements. Plots for two specific functions are shown in \Cref{fig:plots}. There, we see the discrete points accumulate in the less-smooth regions of the functions.

%----------------------------------------------------------------------------------------
%----------------------------------------------------------------------------------------

\subsection{Experiments -- Spline Comparison}

The primary advantage of our algorithm is its nonlinear nature, adapting to the smoothness of a function in each  region. Here, we present a comparison with uniform B-splines. We consider the function
\begin{equation}\label{eq:sine}
    f_\delta(x)
    =
    \frac{1}{x+\delta}
\end{equation}
on the domain $[0,1]$ for small values of $\delta$. We approximate this function using graphs and cubic B-splines, evaluating how many knots are needed to establish an average $L_2$ error below $10^{-3}$.

We see that as $\delta$ decreases, the nonsmooth region near the origin is captured with fewer function evaluations using the adaptive graph spline algorithm.

\begin{table}[ht]
\centering
\begin{tabular}{c || cccc | }
& \multicolumn{4}{c}{Graph Spline}   \\
$\delta$  & num\_init & nfev & refs & $L_2$ error \\ \hline \hline
0.2 & 19 & 89 & 8 & 1.396e-4\\
0.1 & 19 & 137 & 8 & 2.651e-4\\
0.05 & 19 & 211 & 8 & 2.924e-4\\
0.01 & 29 & 489 & 8 & 9.669e-4 \\
\end{tabular}
\begin{tabular}{ |cc  }
\multicolumn{2}{c}{Cubic B-Spline}   \\
knots & $L_2$ error \\ \hline \hline
60 & 8.890e-4\\
60 & 9.539e-4\\
490 & 9.746e-4\\
5910 & 9.917e-4\\
\end{tabular}
\caption{ \textbf{B-Spline Comparison} 
B-spline
}
\label{table:b_spline}
\end{table}

\section{Optimization}
\label{sec:optimization}
\subsection{Algorithm}
Here, we present one potential application of the
discretization described above. We use the discrete approximation to find local minima of a function. The algorithm is straightforward. We create a fine-grid graph representation of the domain. Then, we use the known function values to create an approximation over the entire domain using interpolation. From this approximation, we identify points that are locally minimal.

\subsection{Experiments -- Optimization}

Our experiment is performed on the same set of functions from the previous section. In fact, we use the same discrete approximations to start from. For each function, we build the fine-grid approximation and report the number of local minima found,  the location of the true global minimum, and the location of the nearest local minimum in the approximation. These are reported in \Cref{tab:min_results_generic}.

We see, in each case, that there is a local minimum near the global minimum. 
Using, these initial results, we can find the global minimum to higher precision. Since we are in a one-dimensional setting and the global min is close to a local min, we restrict to a neighborhood of the local minima and apply Brent's method to find the global minimum \cite{brent1973,recipes2007}. For the functions considered, we compute the global minima to an accuracy of $10^{-6}$, using at most 10 additional function evaluations.

\begin{table}[ht!]
\centering
\begin{tabular}{c||c|c|c}
index& number local & true global & nearest local     \\ \hline \hline
2 & 3 & 5.146 & 5.156  \\ \hline 
3 & 20 & -6.775 & -6.787  \\ \hline 
4 & 1 & 2.868 & 2.854  \\ \hline 
5 & 6 & 0.966 & 0.964  \\ \hline 
6 & 14 & 0.680 & 0.703  \\ \hline 
7 & 3 & 5.200 & 5.189  \\ \hline 
8 & 20 & -7.084 & -7.090  \\ \hline 
9 & 3 & 17.039 & 16.987  \\ \hline 
10 & 2 & 7.979 & 7.930  \\ \hline 
11 & 2 & 4.159 & 4.159  \\ \hline 
12 & 3 & 3.142 & 3.160  \\ \hline 
13 & 1 & 0.707 & 0.708  \\ \hline 
14 & 6 & 0.225 & 0.223  \\ \hline 
15 & 2 & 2.414 & 2.368  \\ \hline 
18 & 1 & 2.000 & 2.013  \\ \hline 
20 & 16 & 1.195 & 1.182  \\ \hline 
21 & 6 & 4.795 & 4.814  \\ \hline 
22 & 10 & 14.137 & 14.160  \\ \hline 
\end{tabular}
\caption{Number of local minima identified and comparison with global minimum. The locations of the minima are reported.}
\label{tab:min_results_generic}
\end{table}

\section{Summary and Future Work}
\label{sec:summary}

Here, our focus has been entirely one-dimensional to clearly illustrate the algorithm and show its potential.  We see in the experiments that it is able to adaptively refine a discrete approximation to a function and produce a fine-grid approximation on the full domain, using a limited number of function values.

A key point is that the algorithm inherently computes more function values where a function is less smooth, including critical points of the function, making it useful for optimization problems.

The graph structure of our approximation is not limited to purely one-dimensional objects. In fact, an analogous algorithm can be defined, almost directly, for $n$-dimensional functions. Furthermore, the same key principles can be applied to discretizing functions on smooth manifolds, with suitable modifications on refinement. We leave these and additional extensions for a future work.

\bibliographystyle{IEEEtran}
\bibliography{graph.bib}

\end{document}